\documentclass[11pt]{article}
\usepackage{latexsym}
\usepackage{amsfonts}
\usepackage{mathrsfs,amsthm,amsmath}
\usepackage{color}
\usepackage{todonotes}
\newtheorem{theorem}{Theorem}[section]
\newtheorem{proposition}[theorem]{Proposition}
\newtheorem{definition}[theorem]{Definition}
\newtheorem{corollary}[theorem]{Corollary}
\newtheorem{lemma}{Lemma}[section]

\newtheorem{remark}{Remark}[section]
\newtheorem{example}{Example}[section]

\def\cl#1{{\mathscr #1}}

\def\P{{\mathbb P}}
\def\R{{\mathbb R}}

\def\N{{\mathbb  N}}

\def\E{{\mathbb E}}

\def\Cov{{\rm Cov }}

\def\<{\langle}
\def\>{\rangle}

\begin{document}
\title{Large deviations for perturbed Gaussian processes and 
logarithmic asymptotic estimates for some exit probabilities}
\author{Claudio Macci\thanks{Address: Dipartimento di Matematica,
Universit\`a di Roma Tor Vergata, Via della Ricerca Scientifica,
I-00133 Rome, Italy. e-mail: \texttt{macci@mat.uniroma2.it}}\and
Barbara Pacchiarotti\thanks{Address: Dipartimento di Matematica,
Universit\`a di Roma Tor Vergata, Via della Ricerca Scientifica,
I-00133 Rome, Italy. e-mail: \texttt{pacchiar@mat.uniroma2.it}}}
\maketitle
\begin{abstract}
The main results in this paper concern large deviations for families of non-Gaussian processes obtained as suitable 
perturbations of continuous centered multivariate Gaussian processes which satisfy a large deviation principle. We 
present some corollaries and, as a consequence, we obtain logarithmic asymptotic estimates for exit probabilities 
from suitable halfspaces and quadrants.
\end{abstract}

\section{Introduction}
The theory of large deviations provides a collection of techniques that allow to compute small probabilities of rare events on an
exponential scale (see e.g. \cite{DemboZeitouni} as a reference of this topic). An important feature of large deviations of 
(continuous) Gaussian processes is that the rate functions are expressed in terms of norms of suitable reproducing kernel 
Hilbert spaces associated to their covariance functions. The main results in this paper (Theorems \ref{th:main-equal} and 
\ref{th:main-ind}) concern a general setting and allow to generalize some results in \cite{Pac-Pigl} (see Remark \ref{rem:Pac-Pigl}
for more details). These main results concern families of non-Gaussian processes (on some finite time interval $[0,T]$), obtained 
as suitable perturbations of continuous centered $p$-variate Gaussian processes $(X^n)_{n \in \N}$ which satisfy a large deviation 
principle as $n\to\infty$ (the literature on perturbations of stochastic processes is wide; here we recall 
\cite{DuLeonenkoMaShu} in which the authors consider a perturbation of a Gaussian random field which has some analogy with the one 
in this paper). In particular, in order to have results with more explicit rate function expressions, we assume that the
processes $(X^n)_{n \in \N}$ have independent components (see Condition (C1)). Our approach is motivated by potential connections 
with some processes in the literature and, to this aim, we present some consequences of the main results (Corollaries 
\ref{cor:main-equal} and \ref{cor:main-ind}).

Among the processes in the literature that can be related to our results we cite the \emph{generalized grey Brownian motion}, which 
is a well-known example of model with random diffusivity (see e.g. \cite{SposiniChechkinSenoPagniniMetzler}). Here for completeness 
we recall some other references in the literature on generalized grey Brownian motion. We start with \cite{MuraMainardi} and 
\cite{MuraPagnini} which concern the univariate case (i.e. the case $p=1$ in this paper). Two references on the 
multivariate case are \cite{BockDesmettreDasilva} and \cite{BockGrothausOrge}; some connections between the univariate case and the 
multivariate case are highlighted in \cite{BockDesmettreDasilva} (see eqs. (8) and (9) in that reference; see also \cite{MentrelliPagnini}
cited therein).

In the final part of the paper we study the asymptotic behavior of some exit probabilities. We mean logarithmic asymptotic estimates
(in the fashion of large deviations) of exit probabilities from halfspaces (see Section \ref{sect:LCP-halfspace}) and quadrants (see
Section \ref{sect:LCP-quadrant}). In both cases we obtain the limit in (\ref{eq:LCP-limit}), and we compute its value by using the 
Lagrange multipliers method (this is a standard method used to prove similar results for Gaussian processes in the literature; see for 
instance the results for univariate processes in \cite{Bal-Pac}). Moreover, again in both cases, the exit probabilities can be 
interpreted as level crossing probabilities over a finite time horizon; we recall that level crossing probabilities, or equivalently 
the distributions of first passage times, for univariate generalized grey Brownian motions are studied in the literature (see e.g. 
\cite{SposiniChechkinMetzler}; moreover see \cite{Metzler} for a more general discussion concerning also other topics and models). 
We also recall that, when we deal with the exit probabilities from halfspaces, we have linear combinations of the components of the 
$p$-variate processes, and therefore in some sense we deal with univariate processes. Moreover, when we deal with the exit probabilities
from quadrants, we mean that all the components of the multivariate process reach or exceed (even not simultaneously) 
certain levels $x_1,\ldots,x_q>0$ in a time interval $[0,T]$.

The occurrence of the \emph{exit} event has a natural interpretation in risk theory; for instance, for insurance models, the \emph{ruin}
occurs when a \emph{level crossing} happens (and a level crossing event can be interpreted as an exit event). For instance (here for simplicity
we refer to the exit probabilities from quadrants in this paper) one can consider
an insurance model with $p$ lines of business on some time interval $[0,T]$, and the ruin occurs if we have
$$A_i(t)-c_it\geq x_i\quad (\mbox{for some $i\in\{1,\ldots,p\}$ and $t\in[0,T]$}),$$
where $A_1,\ldots,A_p$ are (possibly correlated) aggregate claim processes, $c_1,\ldots,c_p>0$ are constant premium rates,
and $x_1,\ldots,x_p>0$ are the initial capitals for each line of business. Typically the aggregate claim processes are suitable compound 
sums with positive summands (which represent the claim sizes); however the compound sums could be replaced with some
diffusion approximations, and therefore one has continuous processes (often Brownian motions with some correlations). So, as happens
in several references in the literature, in this paper we are interested in logarithmic asymptotic estimates for sequences of exit 
probabilities $(p_n)_{n\in\N}$ as the limit in (\ref{eq:LCP-limit}) in this paper. One can find this kind of estimates in 
\cite{DebickiJiRolskiRISKS} and in \cite{DebickiKosinskiMandjesRolskiSPA}, where the levels to cross go to infinity with $n$; on the
contrary in this paper the levels are fixed. The literature on level crossing probabilities is wide; here we recall \cite{DebickiHashorvaMichnaJAP}
as a reference with exact (i.e. non-logarithmic) asymptotic estimates, \cite{DebickiHashorvaKrysteckiSAJ} as a reference with approximations 
and bounds, and \cite{BisewskiDebickiKriukov} as a reference with tight bounds and exact asymptotics.

We conclude with the outline of the paper. We start with some preliminaries in Section \ref{sect:prel}. In Section \ref{sect:main} we
present the main results and their corollaries. In Section \ref{sect:LCP} we compute logarithmic asymptotic estimates of exit probabilities
from suitable halfspaces and quadrants. Some inequalities between large deviation rate functions and between logarithmic asymptotic estimates
are obtained in Section \ref{sect:final}, and we discuss some possible connections with the inequalities provided by the supermodular order. 
Finally, in Section \ref{sect:discussion-C2}, we present a discussion on Condition (C2) and how our results can be related to the generalized
grey Brownian motion.

\section{Preliminaries}\label{sect:prel}
In this section we briefly recall some main facts related to reproducing kernel Hilbert space (RKHS) and large deviations for 
Gaussian measures on the Banach space of the continuous functions. 
For a detailed development of this very wide theory we can refer, for example, to the following classical references:
for large deviations, see Section 3.4 in \cite{DeuschelStroock} and Chapter 4 (in particular Sections 4.1, 4.2 and 4.5) in 
\cite{DemboZeitouni}; for reproducing kernel Hilbert spaces, see Chapter 4 (in particular Section 4.3) in \cite{Hid-Hit},
Chapter 2 (in particular Sections 2.2 and 2.3) in \cite{Ber-Tho}, Section 2 in \cite{Mic-Gla}, Section 3 in \cite{Alv-Ros-Law} 
and Section 4.3 in \cite{Lifshits}.

We introduce the setting that we are going to consider throughout the paper.
From now on, given $ T>0 $ and $p\geq 1$, we will denote with $ C([0,T],\R^p) $  the  space of $\R^p$-valued continuous functions 
on $ [0,T] $ endowed with the topology induced by the sup-norm $\|\cdot\|_\infty$, i.e for $f=(f_1,\ldots,f_p)$ then
$$\|f\|_\infty=\sup_{0\leq t \leq T}||f(t)||,$$
where $||\cdot||$ is the euclidean norm in $\R^p$.
The dual set of $ C([0,T],\R^p) $ is the set of vector measures $\lambda=(\lambda_1,\ldots,\lambda_p)$  on $ [0,T] $ and we will 
denote it with $ \mathscr{M}^p[0,T] $. The action of $ \mathscr{M}^p[0,T] $ on $ C([0,T],\R^p) $ is denoted by
$$ \langle\lambda,h\rangle=  \sum_{i=1}^p\int_0^Th_i(t)\, d\lambda_i(t) $$
for every $\lambda\in \cl{M}^p[0,T]$ and $ h=(h_1,\ldots,h_p)\in C([0,T],\R^p)$.
In what follows, we will always suppose our processes to be continuous.

A Gaussian process  $(U(t))_{t \in [0,T]}$  is characterized by the mean function and the covariance function, i.e.
$$m:[0,T]\to \R^p, \quad m_i(t)=\E[U_i(t)]\quad i=1,\ldots,p,$$ and
$$k: [0,T]\times[0,T]\to \R^{p\times p},\quad k_{ij}(t,s)=\Cov(U_i(t),U_j(s))\quad i,j=1,\ldots,p.$$
Recall that the covariance function $ k $ of any Gaussian process is a symmetric, positive definite  function. Since there is a 
one-to-one correspondence between centered Gaussian processes and their covariance functions, we can talk of RKHS relative to a 
positive definite kernel $ k $.
Then let $(U(t))_{t \in [0,T]}$ be a continuous centered Gaussian process, with covariance function $ k $ as above, 
defined on a probability space $(\Omega, \mathscr{F}, \mathbb{P})$. Moreover, define the set
$$ \mathscr{D}=\left\{x \in C([0,T],\R^p): \,x(t)=\int_0^Tk(t,s)\, d\lambda(s), \, \lambda \in  \cl{M}^p[0,T]\right\},$$
where $x(t)=\int_0^T k(t,s)d\lambda(s)$ means
$$x_i(t)=\sum_{j=1}^p\int_0^T k_{ij}(t,s)d\lambda_j(s)\ (\mbox{for all}\ i=1,\ldots,p).$$
As we shall see (see Remark \ref{rem:new} just after Definition \ref{def:RKHS}) the RKHS $\cl H$ relative to the kernel $ k $ can be constructed 
as the completion of the set $ \mathscr{D} $ with respect to a suitable norm. Consider the set of (real) Gaussian random variables
$$ \Gamma=\left\{Y: \, Y=\langle\lambda, U\rangle,\, \lambda \in  \cl{M}^p[0,T]\right\}\subset \mathrm{L}^2(\Omega, \mathscr{F}, \mathbb{P}). $$
We have that, for $ Y_1,Y_2 \in \Gamma, $ say $ Y_i=\langle\lambda^i,U\rangle, $ $ i=1,2, $
$$\begin{array}{ccl}
	\langle Y_1,Y_2\rangle_{\mathrm{L}^2(\Omega, \mathscr{F}, \mathbb{P})}&=& \mbox{Cov}(Y_1,Y_2)\\
	&=& \mbox{Cov}\left( \displaystyle\int_0^T\sum_{j=1}^pU_j(t)\, d\lambda^1_j(t), \displaystyle\int_0^T \sum_{\ell=1}^pU_\ell(t)\, d\lambda^2_\ell(t)\right)\\
	&=& \displaystyle\int_0^T\int_0^T \sum_{j,\ell=1}^pk_{j\ell}(t,s)d\lambda^1_j(t)d\lambda^2_\ell(s).
\end{array}$$
Moreover, define the set
$$  H=\overline{\Gamma}^{\lVert \cdot\rVert_{\mathrm{L}^2(\Omega, \mathscr{F}, \mathbb{P})}} .$$
Then, since $ \mathrm{L}^2 $-limits of Gaussian random variables are still Gaussian, we have that $  H $ is a closed subspace of 
$\mathrm{L}^2(\Omega, \mathscr{F}, \mathbb{P})  $ consisting of real Gaussian random variables. Moreover, it becomes a Hilbert space when endowed 
with the inner product
$$ \langle Y_1,Y_2\rangle_{ H}:= \langle Y_1,Y_2\rangle_{\mathrm{L}^2(\Omega, \mathscr{F}, \mathbb{P})} \quad Y_1,Y_2 \in  H. $$
Consider now the following map,

\begin{eqnarray*}
	\mathscr{S}&:&H \longrightarrow C([0,T], \R^p)\\
	\nonumber&& Y \longmapsto (\mathscr{S}Y)({\cdot})=\mathbb{E}(U({\cdot})Y).
\end{eqnarray*}

The map $\cl S$ is the Lo\'eve isometry (see Theorem 35 in \cite{Ber-Tho}); this isometry has some analogy with the well-known Ito isometry.

\begin{definition}\label{def:RKHS}
	Let $ U=(U(t))_{t \in [0,T]} $ be a continuous centered Gaussian process. We define the \textbf{reproducing kernel Hilbert space} relative to the Gaussian process $ U $ as
	$$ \mathscr{H}:=\mathscr{S}(H)=\{h \in C([0,T],\R^p): \,\,h(t)=(\mathscr{S}Y)(t), Y \in H\} $$
	with an inner product defined as
	$$ \langle h_1,h_2\rangle_{\mathscr{H}}:= \langle \mathscr{S}^{-1}h_1,\mathscr{S}^{-1}h_2\rangle_H=\langle \mathscr{S}^{-1}h_1,\mathscr{S}^{-1}h_2\rangle_{\mathrm{L}^2(\Omega, \mathscr{F}, \mathbb{P})}, \quad h_1,h_2 \in \mathscr{H}. $$
\end{definition}

Now we present some remarks and an example.

\begin{remark}\label{rem:new} \rm
	For any $ \lambda \in \mathscr{M}^p[0,T], $ and $ t\in [0,T] $
	$$ \left(\mathscr{S} \langle\lambda, U\rangle \right)(t)= 
	\mathbb{E}\left( U(t)\int_0^T\sum_{j=1}^pU_j(s)d\lambda_j(s)\right)=\int_0^Tk(t,s)d\lambda(s) $$
	and thus $ \mathscr{S}(\Gamma)=\mathscr{D}$. Then, since $ \Gamma $ is dense in $ H $, and $ \mathscr{S} $ is an isometry,
	we have that
	$$ \mathscr{H}=\mathscr{S}(H)=\overline{\mathscr{S}(\Gamma)}^{\lVert \cdot\rVert_{\mathscr{H}}}=\overline{\mathscr{D}}^{\lVert \cdot\rVert_{\mathscr{H}}}. $$
\end{remark}

\medskip 
\begin{remark}\label{rem:independent-components}\rm
	If the components of the process $U$ are independent, then 
	the set  $\cl D$ is
	\begin{equation}\label{eq:dense-set}
		\mathscr{D}=\left\{x \in C([0,T],\R^p): \,x_i(t)=\int_0^Tk_{ii}(t,s)\, d\lambda_i(s), \,i=1,\ldots, p, \, \lambda\in  \cl{M}^p[0,T]\right\}
	\end{equation}
	and, for $x\in \cl  D$, we have 
	$$
	||x||_{\cl H}^2= \int_0^T\int_0^T \sum_{i=1}^pk_{ii}(t,s)d\lambda_i(t)d\lambda_i(s).
	$$
\end{remark}

Now we present some basic definitions on large deviations.

\begin{definition} \label{def:ldp}
	Let $E$ be a topological space,   $\cl{B}(E)$ the Borel $\sigma$-algebra and  $(\mu_n)_{n\in \N}$ a family of probability measures on $\cl{B}(E)$; let $\gamma \, : \N \rightarrow \mathbb{R}^+ $ be a function, such that  $\gamma_n \rightarrow +\infty$ as $n\to +\infty$. A lower semicontinuous function $I:E\rightarrow [0,+\infty]$ is 
	called rate function. Then we say that $(\mu_n)_{n\in \N}$ satisfies a large deviation principle (LDP) on $E$ with the rate function $I$ and the speed $\gamma_n$ if
	$$
	-\inf_{x \in {\Theta} } I(x) \le \liminf_{n\to+\infty} \frac 1{\gamma_n} \log \mu_n (\Theta)
	$$ 
	for any open set $\Theta$, and
	\begin{equation} \label{eq:upperbound} \limsup_{n\to +\infty}\frac 1 {\gamma_n} \log \mu_n (\Gamma) \le -\inf_{x \in {\Gamma}} I(x)
	\end{equation}
	for any closed set $\Gamma$.
\end{definition}

A rate function $I$ is said \emph{good} if the sets $\{I\le a\}$ are compact for all $a \ge 0$. We also recall the statement of
the contraction principle (see e.g. Theorem 4.2.1 in \cite{DemboZeitouni}).

\begin{proposition}\label{prop:contraction-principle}
        Let $\mathcal{X}$ and $\mathcal{Y}$ be Hausdorff topological spaces 
		and $f:\mathcal{X}\to\mathcal{Y}$ a continuous function. Consider a good rate 
		function $I:\mathcal{X}\to[0,\infty]$.
		\begin{itemize}
			\item For each $y\in\mathcal{Y}$, define
			$$I^\prime(y):=\inf\{I(x) : x\in\mathcal{X}, y= f(x)\}.$$
			Then $I^\prime$ is a good rate function on $\mathcal{Y}$.
			\item If $I$ controls the LDP associated with a family of probability measures
			$\{\mu_n\}$ on $\mathcal{X}$, then $I^\prime$ controls the LDP associated with
			the family of probability measures $\{\mu_n\circ f^{-1}\}$ on $\mathcal{Y}$.
		\end{itemize}
\end{proposition}

In this paper we present large deviation principles for perturbed continuous and centered Gaussian processes. In order
to do that we always suppose that the following Condition (C1) holds; more precisely we mean that an LDP holds and, keeping in mind the principal 
results on large deviations for Gaussian measures (see e.g. \cite{DeuschelStroock}), it is governed by a suitable \emph{quandratic} rate function.

\paragraph{Condition (C1).}
$((X^n(t))_{t \in [0,T]})_{n \in \N} $ is a family of $p$-dimensional, continuous centered Gaussian processes, starting from zero, 
with independent components, which satisfies an LDP with speed $\gamma_n$ and good rate function  
\begin{equation*}
	I(h)=\sum_{i=1}^p I_i(h_i).
\end{equation*}
where $I_i$ is defined by
$$I_i(h_i)=\begin{cases} \frac{1}{2} \left \| h_i \right \|^2_{{\cl  H}_i} & h_i \in{\cl  H_i }\\ +\infty & h_i \notin \cl{H_i}, \end{cases}
\quad (\mbox{for}\ h_i\in C([0,T],\R)),$$
\noindent
where $ \cl{H_i}\subset C([0,T],\R)$ is the RKHS associated to some covariance function $k_{ii}$.

\medskip
\noindent
We can simply write $I(h)=\frac 12 \|h\|^2_{{\cl H}}$ for a 
function $h\in C([0,T], \R^p)$ if we take (as a slight abuse of notation) $\|h\|^2_{{\cl H}}=+\infty$
when $h\notin {\cl H}$.


\section{Main results and corollaries} \label{sect:main}
In this section we present two main results (Theorems \ref{th:main-equal} and \ref{th:main-ind}), and some consequences
(Corollaries \ref{cor:main-equal} and \ref{cor:main-ind}). In particular the main results have some analogies
with the results in \cite{Pac-Pigl}; this fact will be briefly discussed in Remark \ref{rem:Pac-Pigl}, just after 
Theorem \ref{th:main-ind}.

We start with the main results in which we consider a sequence $((A^n,B^n))_{n \in \N}$ such that: it is independent of the 
sequence $((X^n(t))_{t \in [0,T]})_{n \in \N}$ in Condition (C1), it satisfies a LDP with good rate function $I_{(A,B)}$ and 
the speed $\gamma_n$, and $(B^n)_{n \in \N}$ have paths in $C([0,T], \R^p)$. As far as the sequence $(A_n)_{n \in \N}$ is
concerned, we have two cases:
\begin{itemize}
	\item the sequence $(A^n)_{n \in \N}$ in Theorem \ref{th:main-equal} is $[0,+\infty)$ valued;
	\item the sequence $(A^n)_{n \in \N}$ in Theorem \ref{th:main-ind} is $[0,+\infty)^p$ valued.
\end{itemize}
Furthermore, in view of Theorem \ref{th:main-ind}, for $u=(u_1,\ldots,u_p)$ and $v=(v_1,\ldots,v_p)$ in $\R^p$ (from now on we 
will use this notation for possibly random vectors in $\R^p$), we use the symbol $u\circ v$ for the Hadamard product 
between $u$ and $v$, i.e.
$$u\circ v= (u_1 v_1,\ldots,u_pv_p).$$

\begin{theorem}\label{th:main-equal}
	Let $((X^n(t))_{t \in [0,T]})_{n \in \N} $ be as in Condition (C1). Moreover, we consider the family $((A^n,B^n),Z^n)_{n \in \N}$, 
	where  $((A^n,B^n))_{n \in \N}$ is a $[0,+\infty)\times C([0,T], \R^p)$-valued random sequence (thus $(A^n)_{n\in \N}$
	is a sequence of random variables taking values in $[0,+\infty)$, and $(B^n)_{n\in \N}$ is a sequence of continuous stochastic process
	taking values in $\R^p$), independent of $(X^n)_{n\in \N}$. Then we set
	$$Z^n=A^n X^n +B^n\quad (n\in\N),$$
	i.e. $Z^n(t)=A^n X^n(t) +B^n(t)$ for $t\in[0,T]$. We also assume that  $((A^n,B^n))_{n \in \N}$ satisfies a LDP  with the good 
	rate function $I_{(A,B)}$ and the speed $\gamma_n$. Then $((A^n,B^n),Z^n)_{n \in \N}$ satisfies a LDP with the speed $\gamma_n$ 
	and the rate function
	$$I_{(A,B),Z}((a,b),z)=I_{(A,B)}(a,b) + \sum _{i=1}^p J_i(z_i|(a,b_i)),$$
	where 
	\begin{equation}\label{eq:J-i-2}
		J_i(z_i|(a,b_i))= \left\lbrace \begin{aligned}
			&\frac{1}{2a^2}\lVert z_i-b_i\rVert^{2}_{{\cl{H}_{i}}} & 
			a>0\\
			&0&  z_i-b_i=0 \,\mbox{and } a=0\\
			&\infty  & otherwise. 
		\end{aligned} \right.
	\end{equation}
	Moreover $(Z^n)_{n \in \N}$ satisfies a LDP with the speed  $\gamma_n$  and the good rate function
	$$I_Z(z)= \inf_{(a,b) \in[0,+\infty)\times C([0,T], \R^p) }\left\{I_{(A,B)}(a,b) +\sum _{i=1}^p J_i(z_i|(a,b_i))\right\}.$$
\end{theorem}
\begin{proof}
	By the hypotheses we can say that the sequence $((A^n,B^n),X^n)_{n \in \N}$ satisfies a LDP with the speed $\gamma_n$ and 
	the rate function $I_{(A,B),X}$ defined by
	$$I_{(A,B),X}((a,b),x):=I_{(A,B)}(a,b)+I(x),$$
	where $I$ is the rate function in Condition (C1). Then the desired LDPs hold as a consequence of the contraction principle
	(Proposition \ref{prop:contraction-principle}) because the function $(a,b,x)\mapsto ax+b$ defined on
	$[0,+\infty)\times C([0,T], \R^p)\times C([0,T], \R^p)$	is continuous.
\end{proof}

\begin{theorem}\label{th:main-ind}
Let $((X^n(t))_{t \in [0,T]})_{n \in \N} $ be as in Condition (C1). Moreover, we consider the family $((A^n,B^n),Z^n)_{n \in \N}$, 
where  $((A^n,B^n))_{n \in \N}$ is a $[0,+\infty)^p\times C([0,T], \R^p)$-valued random sequence (thus $(A^n)_{n\in \N}$
is a sequence of random variables taking values in $[0,+\infty)^p$, and $(B^n)_{n\in \N}$ is a sequence of continuous stochastic process
taking values in $\R^p$), independent of $(X^n)_{n\in \N}$. Then we set
$$Z^n=A^n \circ X^n +B^n\quad (n\in\N),$$
i.e. $Z^n(t)=A^n\circ X^n(t) +B^n(t)$ for $t\in[0,T]$. We also assume that  $((A^n,B^n))_{n \in \N}$ satisfies a LDP  with the 
good rate function $I_{(A,B)}$ and the speed $\gamma_n$. Then $((A^n,B^n),Z^n)_{n \in \N}$ satisfies a LDP with the speed 
$\gamma_n$ and the rate function $I_{(A,B),Z}$ defined by
$$I_{(A,B),Z}((a,b),z)= I_{(A,B)}(a,b) + \sum _{i=1}^p J_i(z_i|(a_i,b_i)),$$
where
\begin{equation}\label{eq:J-i}
	J_i(z_i|(a_i,b_i))= \left\lbrace \begin{aligned}
		&\frac{1}{2a_i^2}\lVert z_i-b_i\rVert^{2}_{{\cl{H}_{i}}} &  
		a_i>0\\
		&0&  z_i-b_i=0 \,\mbox{and } a_i=0,\\ 
		&\infty  & otherwise.
	\end{aligned} \right.
\end{equation}
\noindent
Moreover $(Z^n)_{n \in \N}$ satisfies a LDP with the speed $\gamma_n$ and the good rate function
$$I_Z(z)= \inf_{(a,b) \in[0,+\infty)^p\times C([0,T], \R^p) }\left\{I_{(A,B)}(a,b) + \sum _{i=1}^p J_i(z_i|(a_i,b_i))\right\}.$$
\end{theorem}
\begin{proof}
	The proof is similar to the proof of Theorem \ref{th:main-equal}. We omit the details.
\end{proof}

\begin{remark}\label{rem:Pac-Pigl}\rm 
The results in this section have some relationship with the results in \cite{Pac-Pigl}.
In that reference it is proved a LDP for a family $((A^n,B^n),Z^n)_{n \in \N}$, where $Z_n=A_nX_n+B_n$ (as we see
in this case we can avoid to refer to the Hadamard product), and the following hypotheses hold:
\begin{itemize}
	\item $((X^n(t))_{t \in [0,T]})_{n \in \N}$ is a family of continuous univariate processes as in Condition (C1) (and
	therefore we can neglect the hypothesis of independent components);
	\item $(A^n,B^n)_{n \in \N}$ is a family of continuous processes with paths in $ C_\alpha([0,1],\R)\times C([0,1],\R)$ 
	(where $C_\alpha([0,1],\R) =\{y\in{C}([0,1],\R) : y(t)\geq \alpha, \,\, t\in[0,1]\}$, equipped with the uniform norm on
	compact sets) which satisfies the same hypotheses as in Theorem \ref{th:main-ind}.
\end{itemize}
So, in order to explain the differences with the model studied in \cite{Pac-Pigl}, we can say that in this paper 
$(A^n)_{n \in \N}$ could be equal to zero, or arbitrarily close to zero, but it is a sequence of random
variables (and not a family of stochastic processes); moreover in this paper $(A^n)_{n \in \N}$ and $((X^n(t))_{t \in [0,T]})_{n \in \N}$
can be multivariate (and not univariate).
\end{remark}

Now we present Corollaries \ref{cor:main-equal} and \ref{cor:main-ind}. To prepare for these corollaries we present the following 
Condition (C2).

\paragraph{Condition (C2).}
$(\widetilde{A}^n)_{n \in \N}$ is a family of positive random variables that satisfies the LDP on $[0,+\infty)$ with the 
speed $\gamma_n$, and good rate function ${\cal J}$ defined by
\begin{equation}\label{eq:rf-C2}
	{\cal J}(x):=dx^\theta\quad (\mbox{for}\ x\geq 0)
\end{equation}
for some $d,\theta>0$.\\
\ \\
This condition will be discussed in Section \ref{sect:discussion-C2}. Moreover, for $d,\theta>0$ in Condition (C2), we consider 
the notation
$$\psi(d,\theta):=2^{\theta/(\theta+2)}\left(d(d\theta)^{-\theta/(\theta+2)}+\frac{1}{2}(d\theta)^{2/(\theta+2)}\right).$$

Now we present two corollaries that are consequences of Theorems \ref{th:main-equal} and \ref{th:main-ind}. We start
with the corollary of Theorem \ref{th:main-equal}.

\begin{corollary}\label{cor:main-equal}
	Let $((X^n(t))_{t \in [0,T]})_{n \in \N} $ be as in Condition (C1). Moreover, we consider the family $((A^n,B^n),Z^n)_{n \in \N}$, 
	where  $((A^n,B^n))_{n \in \N}$ is a $[0,+\infty)\times C([0,T], \R^p)$-valued random sequence (thus $(A^n)_{n\in \N}$
	is a sequence of random variables taking values in $[0,+\infty)$, and $(B^n)_{n\in \N}$ is a sequence of continuous stochastic process
	taking values in $\R^p$), independent of $(X^n)_{n\in \N}$. Then we set
	$$Z^n=A^n X^n +B^n\quad (n\in\N),$$
	i.e. $Z^n(t)=A^n X^n(t) +B^n(t)$ for $t\in[0,T]$. 
	We also assume that: $(A^n)_{n\in\N}=(\widetilde{A}^n)_{n \in \N}$, where $(\widetilde{A}^n)_{n \in \N}$ is as in 
	Condition (C2) for some $d,\theta>0$; $(B^n)_{n\in\N}$ satisfies a LDP with the speed $\gamma_n$ and the good rate function $I_B$;
	$(A^n)_{n\in\N}$ and $(B^n)_{n\in\N}$ are independent. Then $((A^n,B^n),Z^n)_{n \in \N}$ satisfies a LDP with the speed 
	$\gamma_n$ and the rate function $I_{(A,B),Z}$ defined by
	$$I_{(A,B),Z}((a,b),z)= I_B(b) + d a^\theta + \sum _{i=1}^p J_i(z_i|(a,b_i)),$$
	where $J_i(z_i|(a,b_i))$ is the function defined by (\ref{eq:J-i-2}); so in particular we have
	$$I_{(A,B),Z}((a,b),z)=\left\{\begin{array}{ll}
		I_B(b)+d a^\theta+\sum_{i=1}^p\frac{1}{2a^2}\lVert z_i-b_i\rVert^{2}_{{\cl{H}_{i}}}&\ \mbox{if}\ a>0\\
		I_B(b) &\ \mbox{if}\ a=0 \,\mbox{and}\, z_i-b_i=0, \, i=1,\ldots,p\\
		\infty&\ \mbox{otherwise}.
	\end{array}\right.$$
	Moreover $(Z^n)_{n \in \N}$ satisfies a LDP with the speed $\gamma_n$ and the good rate function
	$$\displaylines{I_Z(z)= \inf_{(a,b) \in[0,+\infty)\times C([0,T], \R^p) }\left\{I_B(b) + d a^\theta + \sum _{i=1}^p J_i(z_i|(a,b_i))\right\}\cr
		=\inf_{b\in C([0,T], \R^p) }\left\{I_B(b) +  \psi(d,\theta)\left(\sum _{i=1}^p\lVert z_i- b_i\rVert^2_{{\cl{H}_{i}}}\right)^{\theta/(\theta+2)}\right\}.}$$
\end{corollary}
\begin{proof}
	It is an immediate consequence of Theorem \ref{th:main-equal} and some applications of the contraction principle.
\end{proof}

\begin{remark}\label{rem:main-equal}\rm
	In view of the results presented in Section \ref{sect:LCP} it is useful to refer to the sequence $(Z^n)_{n\in \N}$ 
	in Corollary \ref{cor:main-equal} when $(B^n)_{n\in \N}$
	is a constant deterministic sequence equal to some $\hat{b}=(\hat b_1,\ldots,\hat b_p)\in C([0,T],\mathbb{R}^p)$;
	in such a case we have $I_B(b)=0$ if $b=\hat b$, and infinity otherwise, and therefore we can say that $I_Z(z)=I_Z^{(=)}(z)$, where
	\begin{equation}\label{eq:rate-fun-Z-ggG}
		I_Z^{(=)}(z):=\psi(d,\theta)\left(\frac12 \sum_{i=1}^p\lVert z_i-\hat{b}_i\rVert^{2}_{{\cl{H}_{i}}}\right)^{\theta/(\theta+2)}.
	\end{equation}
\end{remark}

Now we present the corollary of Theorem \ref{th:main-ind}.

\begin{corollary}\label{cor:main-ind}
	Let $((X^n(t))_{t \in [0,T]})_{n \in \N} $ be as in Condition (C1). Moreover, we consider the family $((A^n,B^n),Z^n)_{n \in \N}$, 
	where  $((A^n,B^n))_{n \in \N}$ is a $[0,+\infty)^p\times C([0,T], \R^p)$-valued random sequence (thus $(A^n)_{n\in \N}$
	is a sequence of random variables taking values in $[0,+\infty)^p$, and $(B^n)_{n\in \N}$ is a sequence of continuous stochastic process
	taking values in $\R^p$), independent of $(X^n)_{n\in \N}$. Then we set
	$$Z^n=A^n \circ X^n +B^n\quad (n\in\N),$$
	i.e. $Z^n(t)=A^n\circ X^n(t) +B^n(t)$ for $t\in[0,T]$. 
	We also assume that: $(A^n)_{n\in\N}=((\widetilde{A}^n_1,\ldots,\widetilde{A}^n_p))_{n \in \N}$, where 
	$(\widetilde{A}^n_1)_{n \in \N},\ldots(\widetilde{A}^n_p)_{n \in \N},$ are as in Condition (C2) for some $d_1,\theta_1,\ldots,d_p,\theta_p>0$;
	$(B^n)_{n\in\N}$ satisfies a LDP with the speed $\gamma_n$ and the good rate function $I_B$;
	$(A^n)_{n\in\N}$ and $(B^n)_{n\in\N}$ are independent. Then $((A^n,B^n),Z^n)_{n \in \N}$ satisfies a LDP with the speed 
	$\gamma_n$ and the rate function $I_{(A,B),Z}$ defined by
	$$I_{(A,B),Z}((a,b),z)= I_B(b) + \sum _{i=1}^p \{d_i a_i^{\theta_i} + J_i(z_i|(a_i,b_i))\},$$
	where $J_i(z_i|(a_i,b_i))$ is the function defined by (\ref{eq:J-i}); so in particular, if we set
	$${\cal I}^+(a)=\{1\leq i\leq p: a_i>0\}\quad\mbox{and}\quad {\cal I}^0(a)=\{1\leq i\leq p: a_i=0\},$$
	we have
	$$I_{(A,B),Z}((a,b),z)=\left\{\begin{array}{ll}
			I_B(b)+\sum_{i\in {\cal I}^+(a)}\left\{d_i a_i^{\theta_i}+\frac{1}{2a_i^2}\lVert z_i- b_i\rVert^{2}_{{\cl{H}_{i}}}\right\}&
			\ \mbox{if}\ z_i-b_i=0\ \mbox{for}\ i \in {\cal I}^0(a),\\
			\infty&\ \mbox{otherwise}.
		\end{array}\right.$$
	Moreover $(Z^n)_{n \in \N}$ satisfies a LDP with the speed $\gamma_n$ and the good rate function
	$$\displaylines{I_Z(z)=\inf_{(a,b) \in[0,+\infty)^p\times C([0,T], \R^p) }\left\{I_B(b) + \sum _{i=1}^p \{d_i a_i^{\theta_i} + J_i(z_i|(a_i,b_i))\}\right\}\cr
		=\inf_{b\in C([0,T], \R^p) }\left\{I_B(b) + \sum _{i=1}^p \psi(d_i,\theta_i)\lVert z_i- b_i\rVert^{2\theta_i/(\theta_i+2)}_{{\cl{H}_{i}}}\right\}.}$$
\end{corollary}
\begin{proof}
	It is an immediate consequence of Theorem \ref{th:main-ind} and some applications of the contraction principle.
\end{proof}

\begin{remark}\label{rem:main-ind}\rm
	In view of the results presented in Section \ref{sect:LCP} it is useful to refer to the sequence $(Z^n)_{n\in \N}$ 
	in Corollary \ref{cor:main-ind} when $(B^n)_{n\in \N}$
	is a constant deterministic sequence equal to some $\hat{b}=(\hat b_1,\ldots,\hat b_p)\in C([0,T],\mathbb{R}^p)$;
	in such a case we have $I_B(b)=0$ if $b=\hat b$, and infinity otherwise, and therefore we can say that $I_Z(z)=I_Z^{(\mathrm{\bot})}(z)$, 
	where
	\begin{equation}\label{eq:rate-fun-Z}
		I_Z^{(\bot)}(z):=\sum _{i=1}^p \psi(d_i,\theta_i)\left(\frac12 \lVert z_i- \hat b_i\rVert^2_{{\cl{H}_{i}}}\right)^{\theta_i/(\theta_i+2)}.
	\end{equation}
\end{remark}

We can say that we are not aware if Corollaries \ref{cor:main-equal} and \ref{cor:main-ind} (and their consequences in
Remarks \ref{rem:main-equal} and \ref{rem:main-ind}) provide known results. However all those statements are consequences of Theorems
\ref{th:main-equal} and \ref{th:main-ind} which concern more general situations; so all those statements are useful to illustrate 
what happens under the hypotheses that we consider to obtain the results in Sections \ref{sect:LCP} and \ref{sect:final}.

\section{Asymptotic results for some exit probabilities}\label{sect:LCP}
In this section we obtain some asymptotic estimates for two exit probabilities concerning the processes $(Z^n)_{n\in\N}$ presented in 
Corollaries \ref{cor:main-equal} and \ref{cor:main-ind} (together with Remarks \ref{rem:main-equal} and \ref{rem:main-ind}, respectively); 
more precisely, in both cases, we restrict our attention to the case in which $(B^n)_{n\in \N}$ is a constant deterministic sequence (that is 
$B^n=\hat b$ for every $n\in \N$, for some $\hat{b}=(\hat b_1,\ldots,\hat b_p)\in C([0,T],\mathbb{R}^p)$). So the rate function $I_Z$ in 
Corollary \ref{cor:main-equal} coincides with $I_Z^{(=)}$ in (\ref{eq:rate-fun-Z-ggG}), and the rate function $I_Z$ in Corollary 
\ref{cor:main-ind} coincides with $I_Z^{(\bot)}$ in (\ref{eq:rate-fun-Z}).

We consider exit probabilities from halfspaces in Section \ref{sect:LCP-halfspace}, and from quadrants in Section \ref{sect:LCP-quadrant}; 
moreover, in both cases, the \emph{set of paths which lead to exit} will be denoted by $\cl{A}$ and it is a closed set. We deal with two 
sequences of exit probabilities $(p_n)_{n\in\N}$ (see (\ref{eq:LCP-halfspace-def}) for the exit of a halfspace, and (\ref{eq:LCP-quadrant-def})
for the exit from quadrants); then, by taking into account the LDPs in Corollaries \ref{cor:main-equal} and \ref{cor:main-ind}, we have
\begin{equation*}
    \begin{aligned}
			-\inf_{z \in \cl{A}^\circ}I_Z(z)\leq \liminf_{n \rightarrow +\infty}\frac{1}{\gamma_n}\log(p_n)
			\leq \limsup_{n \rightarrow +\infty}\frac{1}{\gamma_n}\log(p_n)\leq -\inf_{z \in \cl{A}}I_Z(z),
		\end{aligned}
\end{equation*}
where $\cl{A}^\circ$ is the interior of the set $\cl{A}$, and $I_Z$ is the rate function $I_Z^{(=)}$ in (\ref{eq:rate-fun-Z-ggG}), 
or the rate function $I_Z^{(\bot)}$ in (\ref{eq:rate-fun-Z}). Moreover, in all the cases studied below, we check that
\begin{equation}\label{eq:coincident-infima}
w:=\inf_{z \in \cl{A}^\circ}I_Z(z)=\inf_{z \in \cl{A}}I_Z(z),
\end{equation}
which yields the logarithmic asymptotic estimate
\begin{equation}\label{eq:LCP-limit}
\lim_{n\to +\infty}\frac{1}{\gamma_n}\log(p_n)=-w.
\end{equation}
The value $w$ is called \emph{exponential decay rate}. In particular we also have $w=I_Z(z^*)$ for some $z^*\in \cl{A}$ and, in the fashion 
of large deviations, $z^*$ is said to be \emph{a most likely path leading to exit} (for some more details on this concept 
see e.g. Lemma 4.2 in \cite{GaneshOConnellWischik}; another reference is \cite{Mandjes}, p. 45). As far as the equality in (\ref{eq:coincident-infima})
is concerned, in general we trivially have
$$\inf_{z \in \cl{A}^\circ}I_Z(z)\geq\inf_{z \in \cl{A}}I_Z(z),$$
and therefore only the inverse inequality has to be checked.

\begin{remark}\label{rem:constraint-on-theta}\rm
	The computation of the exponential decay rate $w$ is a minimization problem concerning the rate function $I_Z$.
	In our results we need a strict convexity property for the rate function $I_Z$, and this allows to say that there exists a unique minimizer.
	In Propositions \ref{prop:LCP-halfspace-equal} and \ref{prop:LCP-quadrant-equal} we refer to $I_Z$ in (\ref{eq:rate-fun-Z-ggG}); then,
	since the basis of the power in (\ref{eq:rate-fun-Z-ggG}) is a strictly convex function, we can choose $d,\theta>0$ without any restriction.
	In Propositions \ref{prop:LCP-halfspace-independent} and \ref{prop:LCP-quadrant-independent} we refer to $I_Z$ in (\ref{eq:rate-fun-Z});
	then we have a strictly convex function if the exponents $\frac{\theta_1}{\theta_1+2},\ldots,\frac{\theta_p}{\theta_p+2}$ are larger than 
	$\frac{1}{2}$ (with $\theta_1,\ldots,\theta_p>0$) or, equivalently, if $\theta_1,\ldots,\theta_p>2$. We also remark that in 
	Proposition \ref{prop:LCP-halfspace-independent} we have $d_1=\cdots=d_p=d>0$ and $\theta_1=\cdots=\theta_p=\theta>2$. Finally we can say that
	the inequality $\theta>2$ comes up in Proposition \ref{prop:inequalities}, where we consider some comparison (inequalities) between rate 
	functions and between exponential decay rates.
\end{remark}

We recall that the processes $(Z^n)_{n\in\N}$ presented in Corollaries \ref{cor:main-equal} and \ref{cor:main-ind} are
defined in terms of a continuous centered Gaussian process $(X(t))_{t\in[0,T]}$ with independent components (see Condition (C1));
so we can refer to the dense set $\cl{D}$ as in (\ref{eq:dense-set}) in Remark \ref{rem:independent-components}. Moreover, in 
what follows we also take into account that
$$\inf_{z \in \cl{A}}I_Z(z)=\inf_{z \in \cl{A}\cap\cl{D}_{\hat b}}I_Z(z),$$
where 
$$\cl{D}_{\hat b}:=\cl{D} + {\hat b}=\{z=y+\hat b: y\in\cl{D}\}.$$
More precisely we can say that $z\in\cl{D}_{\hat b}$ if and only if, for some $\lambda_1,\ldots,\lambda_p\in\cl{M}^1[0,T]$,
we have
\begin{equation}\label{eq:z-dense}
z_i(u) =\int_0^T {k_{ii}}(u,v)\ d\lambda_i(v)+\hat b_i(u), \quad u \in [0,T],\ i=1,\ldots,p;
\end{equation}
then, if $z$ is as in (\ref{eq:z-dense}), we have
\begin{equation}\label{eq:z-dense-norm}
\lVert z_i- \hat b_i\rVert^{2}_{{\cl{H}_{i}}}=\int_0^T \int_0^T k_{ii}(u,v)\ d\lambda_i(u)\ d\lambda_i(v).
\end{equation}

\subsection{Exit probabilities from halfspaces}\label{sect:LCP-halfspace}
In this section we consider
\begin{equation}\label{eq:LCP-halfspace-def}
p_n:=\P\Big(\sup_{t\in[0,T]}\langle Z^n(t),\xi\rangle \geq x\Big),
\end{equation}
for $0\neq\xi\in\R^p$, $\xi_i\geq0$ for $i=1,\ldots,p$, and $x>0$ such that $x-\langle \hat  b(t), \xi \rangle>0$ for every 
$t\in [0,T]$. Note that, in some sense, we deal with a sequence of univariate processes $(\langle Z^n(t),\xi\rangle)_{n\in\N}$. 
Moreover the set $\cl A$ is defined by
$$\cl A:=\left\{ z \in {C}([0,T], \R^p): \sup_{t\in[0,T]}\langle z(t),\xi\rangle \geq x\right\}.$$
Then we have
\begin{equation}\label{eq:union-of-paths-halfspace}
\cl{A}=\bigcup_{t\in[0,T]}\cl{A}_t,\ \mbox{where}\
\cl{A}_t:=\left\lbrace z \in {C}([0,T],\R^p): \langle z(t), \xi \rangle= x\right\rbrace.
\end{equation}

\begin{remark}\label{rem:coincident-infima-halfspace}\rm
Here we check (\ref{eq:coincident-infima}), and therefore we only have to check the inequality
$$\inf_{z \in \cl{A}^\circ}I_Z(z)\leq\inf_{z \in \cl{A}}I_Z(z).$$
There exists $z^*\in \cl{A}$ such that $I_Z(z^*)=\inf_{z \in \cl{A}}I_Z(z)$. If $z^*\in \cl{A}^\circ$ this is
trivial, and therefore here we assume that $z^*\notin \cl{A}^\circ$. In this case we have
$$\sup_{0\leq t \leq T}\langle z^*(t), \xi \rangle= x;$$
so, for some $t^*\in[0,T]$, we have $\langle z^*(t^*), \xi \rangle= x$ and
$$\langle z^*(t^*)-\hat  b(t^*), \xi \rangle= x-\langle \hat  b(t^*), \xi \rangle>0.$$
Then, for every $\varepsilon>0$, let $z^{*,\varepsilon}$ be defined by
$$z^{*,\varepsilon}(t):=(1+\varepsilon)(z^*(t)-\hat  b(t))+\hat  b(t);$$
so $z^{*,\varepsilon}\in\cl{A}^\circ$ because
\begin{multline*}
\sup_{0\leq t \leq T}\langle z^{*,\varepsilon}(t), \xi \rangle
\geq \langle z^{*,\varepsilon}(t^*), \xi \rangle=(1+\varepsilon)\langle z^*(t^*)-\hat  b(t^*), \xi \rangle+\langle \hat  b(t^*), \xi \rangle\\
=(1+\varepsilon)(x-\langle \hat  b(t^*), \xi \rangle)+\langle \hat  b(t^*), \xi \rangle
=x+\varepsilon(x-\langle \hat  b(t^*), \xi \rangle)>x.
\end{multline*}
Finally we get
$$\inf_{z \in \cl{A}^\circ}I_Z(z)\leq I_Z(z^{*,\varepsilon})\to I_Z(z^*)\ \mbox{as}\ \varepsilon\to 0,$$
where the limit can be checked by taking $I_Z$ as $I_Z^{(=)}$ in (\ref{eq:rate-fun-Z-ggG})
or $I_Z^{(\bot)}$ in (\ref{eq:rate-fun-Z}).
\end{remark}

In the proof of the following Propositions \ref{prop:LCP-halfspace-equal} and
\ref{prop:LCP-halfspace-independent} we compute $w$ noting that
\begin{equation}\label{eq:w-def-halfspace}
w=\inf_{z \in \cl{A}}I_Z(z)=\inf_{t\in[0,T]}\inf_{z \in \cl{A}_t}I_Z(z);
\end{equation}
moreover, for every $t \in [0,T]$, we compute $\inf_{z \in \cl{A}_t}I_Z(z)$ by taking $z$ in the dense set $\cl{D}_{\hat b}$
defined above (see eq. (\ref{eq:z-dense})) by applying the Lagrange multipliers methods (we can do this because
the rate function $I_Z$ is strictly convex; see Remark \ref{rem:constraint-on-theta}). Now we are ready to prove Proposition 
\ref{prop:LCP-halfspace-equal} (see Corollary \ref{cor:main-equal} and Remark \ref{rem:main-equal}) and Proposition 
\ref{prop:LCP-halfspace-independent} (see Corollary \ref{cor:main-ind} and Remark \ref{rem:main-ind}; actually we also need the 
further conditions $d_1=\cdots=d_p=d$ and $\theta_1=\cdots=\theta_p=\theta$).

\begin{proposition}\label{prop:LCP-halfspace-equal}
Let $p_n$ be as in (\ref{eq:LCP-halfspace-def}), where
$$Z^n(t)=\widetilde{A}^nX^n(t)+\hat b(t)\quad(t\in[0,T]),$$
$(X^n)_{n\in\N}$ is a continuous Gaussian process as in Condition (C1), $\hat b$ is as in Remark \ref{rem:main-equal}, and 
$(\widetilde{A}^n)_{n\in\N}$ in Condition (C2) for some $d,\theta>0$. Then (\ref{eq:LCP-limit}) holds with $w=w_H^{(=)}$, where
$$w_H^{(=)}:=\psi(d,\theta)
\inf_{t\in[0,T]}\left[\frac 12 \frac{(-\langle \hat  b(t), \xi \rangle+x)^2} {\sum_{j=1}^p \xi_j^{2}k_{jj}(t,t)}\right]^{\theta/(\theta+2)}.$$
\end{proposition}
\begin{proof}
We have to compute $w$ in (\ref{eq:w-def-halfspace}) where $\cl{A}_t$ is as in (\ref{eq:union-of-paths-halfspace}) and
$I_Z$ is the rate function $I_Z^{(=)}$ in (\ref{eq:rate-fun-Z-ggG}). Thus, by referring to the elements in 
the dense set of paths $\cl{D}_{\hat b}$ (see (\ref{eq:z-dense}) and (\ref{eq:z-dense-norm})), we have to minimize
$$\psi(d,\theta)\left(\sum_{i=1}^p\frac{1}{2}\int_0^T \int_0^T k_{ii}(u,v)\ d\lambda_i(u)\ d\lambda_i(v)\right)^{\theta/(\theta+2)}$$
with respect to the vector measures $\lambda$, subjected to the following constraint
$$ \langle z(t), \xi \rangle=\sum_{i=1}^p \xi_i \int_0^T k_{ii}(t,v)\ d\lambda_i(v)+ \langle \hat  b(t), \xi \rangle=x. $$
So we use the method of Lagrange multipliers and, for $\gamma\in\R$, we have to find the stationary points of
$$\displaylines{\cl{L}(\lambda,\gamma):=\psi(d,\theta)\left(\sum_{i=1}^p\frac{1}{2}\int_0^T \int_0^T k_{ii}(u,v)
\ d\lambda_i(u)\ d\lambda_i(v)\right)^{\theta/(\theta+2)}\cr
- \gamma \Big(\sum_{i=1}^p\xi_i\int_0^T k_{ii}(t,v)\  d\lambda_i(v)+\langle \hat  b(t), \xi \rangle-x\Big).}$$
Then, for every  $\eta\in  \cl{M}^1[0, 1]$, we have
$$\displaylines{\psi(d,\theta)\frac{\theta}{\theta+2}\left(\frac{1}{2}\sum_{i=1}^p\int_0^T \int_0^T k_{ii}(u,v)\ d\lambda_i(u)\ d\lambda_i(v)\right)^{-\frac{2}{\theta+2}}\int_0^T\int_0^T k_{ii}(u,v)\ d\lambda_i(u)\ d\eta(v)+
\cr - \gamma \xi_i \int_0^T k_{ii}(t,v) d\eta(v)=0,\quad i=1,\ldots,p;}$$
therefore we have to find a vector measures $\lambda$ such that, for every $v\in [0,T]$,
$$\displaylines{\psi(d,\theta)\frac{\theta}{\theta+2}\left(\frac{1}{2}\sum_{i=1}^p\int_0^T \int_0^T k_{ii}(u,v)
\ d\lambda_i(u)\ d\lambda_i(v)\right)^{-\frac{2}{\theta+2}} \int_0^T k_{ii}(u,v)\ d\lambda_i(u)\cr
- \gamma \xi_i k_{ii}(t,v) =0,\quad i=1,\ldots,p.}$$
A solution is $\lambda_i=c_i\delta_{\{ t\} }$ for some $c_i\in\R$ (for every $i=1,\ldots,p$). 
So we have to consider the following system
$$\left\{
\begin{array}{ll}
\displaystyle \psi(d,\theta)\frac{\theta}{\theta+2}\left(\frac{1}{2}\sum_{i=1}^p c_i^2k_{ii}(t,t)\right)^{-\frac{2}{\theta+2}} c_i-\gamma \xi_i=0& i=1,\ldots,p\\
\displaystyle\sum_{i=1}^p \xi_i c_i k_{ii}(t,t)=-\langle \hat  b(t), \xi \rangle+x,&\end{array}
\right.$$
which can be explicitly solved; in fact we can check that
$$\gamma=\psi(d,\theta)\frac{\theta}{\theta+2}\Big(\frac 12\Big)^{-\frac{2}{\theta+2}}
\frac{(-\langle \hat  b(t), \xi \rangle+x)^{(\theta-2)/(\theta+2)}}{\Big(\sum_{i=1}^p \xi_i^{2}k_{ii}(t,t)\Big)^{\theta/(\theta+2)}}$$
and
$$\displaylines{c_i=  \frac {\xi_i(-\langle \hat  b(t), \xi \rangle+x)}{\sum_{j=1}^p{\xi_j^{2  }k_{jj}(t,t) }},\quad i=1,\ldots,p.}$$
So this is a solution of the Lagrange multipliers problem, and it is therefore a critique point for the functional we want to minimize.
Moreover, since this is a strictly convex good rate function restricted on a linear subspace of $\cl{M}^p[0, 1]$, it is still strictly
convex, and therefore the critique point is actually its unique point of minimum. In conclusion we easily get the desired expression 
for $w$ by (\ref{eq:w-def-halfspace}); indeed we have
$\inf_{z \in \cl{A}_t}I_Z(z)=\psi(d,\theta)\left[\sum_{i=1}^p\frac 12 c_i^2 k_{ii}(t,t)\right]^{\theta/(\theta+2)}$
with $c_1,\ldots,c_p$ computed above.
\end{proof}

\begin{proposition}\label{prop:LCP-halfspace-independent}
Let $p_n$ be as in (\ref{eq:LCP-halfspace-def}), where
$$Z^n(t)=A^n\circ X^n(t)+\hat b(t)\quad(t\in[0,T]),$$
$(X^n)_{n\in\N}$ is a continuous Gaussian process as in Condition (C1), $\hat b$ is as in Remark \ref{rem:main-ind}, and 
$(A^n)_{n\in\N}=((\widetilde{A}^n_1,\ldots,\widetilde{A}^n_p))_{n \in \N}$, where 
$(\widetilde{A}^n_1)_{n \in \N},\ldots(\widetilde{A}^n_p)_{n \in \N},$ are as in Condition (C2) for some $d_1=\cdots=d_p=d>0$ and 
$\theta_1=\cdots=\theta_p=\theta>2$. Then (\ref{eq:LCP-limit}) holds with $w=w_H^{(\bot)}$, where
$$w_H^{(\bot)}:=\psi(d,\theta)
\inf_{t\in[0,T]}\left[\frac 12 \frac{(-\langle \hat  b(t), \xi \rangle+x)^2} {\Big(\sum_{j=1}^p
\xi_j^{2\theta/(\theta-2)}k_{jj}(t,t)^{\theta/(\theta-2)}\Big)^{(\theta-2)/\theta}}\right]^{\theta/(\theta+2)}.$$
\end{proposition}
\begin{proof}
    We remark that the hypotheses $d_1=\cdots=d_p=d$ and $\theta_1=\cdots=\theta_p=\theta$ for some $d>0$ and $\theta>2$ are
    actually useful only in the final part of the proof; so we refer these hypotheses when it will be needed.
    We have to compute $w$ in (\ref{eq:w-def-halfspace}) where $\cl{A}_t$ is as in (\ref{eq:union-of-paths-halfspace}) and
    $I_Z$ is the rate function $I_Z^{(\bot)}$ in (\ref{eq:rate-fun-Z}). We follow the same lines of the proof of 
    Proposition \ref{prop:LCP-halfspace-equal} and we omit some details. Then we use again the method of Lagrange multipliers and,
	for $\gamma\in\R$, we have to find the stationary points of
	$$\displaylines{\cl{L}(\lambda,\gamma):=\sum_{i=1}^p\psi(d_i,\theta_i)\left(\frac{1}{2}\int_0^T \int_0^T k_{ii}(u,v)
		\ d\lambda_i(u)\ d\lambda_i(v)\right)^{\theta_i/(\theta_i+2)}\!\!\!\cr
		- \gamma \Big(\sum_{i=1}^p\xi_i\int_0^T k_{ii}(t,v)\  d\lambda_i(v)+ \langle \hat  b(t), \xi \rangle-x\Big).}$$
	Then, for every  $\eta\in  \cl{M}^1[0, 1]$, we have
	$$\displaylines{\psi(d_i,\theta_i)\frac{\theta_i}{\theta_i+2}\left(\frac{1}{2}\int_0^T \int_0^T k_{ii}(u,v)\ d\lambda_i(u)\ d\lambda_i(v)\right)^{-\frac{2}{\theta_i+2}}\int_0^T\int_0^T k_{ii}(u,v)\ d\lambda_i(u)\ d\eta(v)+
		\cr - \gamma \xi_i \int_0^T k_{ii}(t,v) d\eta(v)=0,\quad i=1,\ldots,p;}$$
	therefore we have to find a vector measures $\lambda$ such that, for every $v\in [0,T]$,
	$$\displaylines{\psi(d_i,\theta_i)\frac{\theta_i}{\theta_i+2}\left(\frac{1}{2}\int_0^T \int_0^T k_{ii}(u,v)
		\ d\lambda_i(u)\ d\lambda_i(v)\right)^{-\frac{2}{\theta_i+2}} \int_0^T k_{ii}(u,v)\ d\lambda_i(u)\cr
		- \gamma \xi_i k_{ii}(t,v) =0,\quad i=1,\ldots,p.}$$
	Then, again, a solution is $\lambda_i=c_i\delta_{\{ t\} }$ for some $c_i\in\R$ (for every $i=1,\ldots,p$). So we have to
	consider the following system
	$$\left\{
	\begin{array}{ll}
		\displaystyle \psi(d_i,\theta_i)\frac{\theta_i}{\theta_i+2}
		\left(\frac{1}{2}k_{ii}(t,t)\right)^{-\frac{2}{\theta_i+2}} c_i^{(\theta_i-2)/(\theta_i+2)}-\gamma \xi_i=0& i=1,\ldots,p\\
		\displaystyle\sum_{i=1}^p \xi_i c_i k_{ii}(t,t)=-\langle \hat  b(t), \xi \rangle+x.&\end{array}
	\right.$$
	From now on we take $d_1=\cdots=d_p=d$ and $\theta_1=\cdots=\theta_p=\theta$ as in the statement of the proposition, and the system can be explicitly
	solved; in fact we can check that
	$$\gamma=\left( \frac {-\langle \hat  b(t), \xi \rangle+x}{\sum_{i=1}^p \xi_i^{2\theta/(\theta-2)}\Big(\frac{1}{\psi(d,\theta)}\Big)^{(\theta+2)/(\theta-2)}
	\Big(\frac 12\Big)^{2/(\theta-2)}k_{ii}(t,t)^{\theta/(\theta-2)}}
	\right)^{(\theta-2)/(\theta+2)}$$
	and
	$$\displaylines{c_i=\gamma^{(\theta+2)/(\theta-2)} \xi_i^{(\theta+2)/(\theta-2)}\Big(\frac{1}{\psi(d,\theta)}\Big)^{(\theta+2)/(\theta-2)}
		\Big(\frac 12k_{ii}(t,t)\Big)^{2/(\theta-2)}\cr
		= (-\langle \hat  b(t), \xi \rangle+x) \frac {\xi_i^{(\theta+2)/(\theta-2)}k_{ii}(t,t)^{2/(\theta-2)}}{\sum_{j=1}^p{\xi_j^{2\theta/(\theta-2)}k_{jj}(t,t)^{\theta/(\theta-2)}} },\quad i=1,\ldots,p.}$$
	So we can conclude following the same lines of the final part of the proof of Proposition \ref{prop:LCP-halfspace-equal}.
	In particular we easily get the desired expression for $w$ by (\ref{eq:w-def-halfspace}); indeed we have 
	$\inf_{z \in \cl{A}_t}I_Z(z)=\psi(d,\theta)\sum_{i=1}^p\left[\frac 12 c_i^2 k_{ii}(t,t)\right]^{\theta/(\theta+2)}$
	with $c_1,\ldots,c_p$ computed above (actually here we need some more computations with respect to the
	case of Proposition \ref{prop:LCP-halfspace-equal}).
\end{proof}

\paragraph{On the most likely paths leading to exit.}
A closer look at the proof of Proposition \ref{prop:LCP-halfspace-equal} reveals that the following function $z^*$ is a most likely path leading to exit:
$$z_i^*(u):=\frac{\xi_i(-\langle \hat  b(t^*), \xi \rangle+x)}{\sum_{j=1}^p{\xi_j^{2 }k_{jj}(t^*,t^*)} } k_{ii}(u,t^*)
+\hat b_i(u), \quad u \in [0,T],\ i=1,\ldots,p,$$
where
$$t^*={\rm argmin}_{t\in[0,T]}\left[\frac 12 \frac{(-\langle \hat  b(t), \xi \rangle+x)^2} {\sum_{j=1}^p \xi_j^{2}k_{jj}(t,t)}\right]^{\theta/(\theta+2)}.$$
Similarly, for Proposition \ref{prop:LCP-halfspace-independent}, we can define $z^*$ as 
follows:
$$z_i^*(u):=(-\langle \hat  b(t^*), \xi \rangle+x)\frac{\xi_i^{2\theta/(\theta-2)-1}k_{ii}(t^*,t^*)^{\theta/(\theta-2)-1}}
{\sum_{j=1}^p{\xi_j^{2\theta/(\theta-2)}k_{jj}(t^*,t^*)^{\theta/(\theta-2)}} } k_{ii}(u,t^*)+\hat b_i(u), \quad u \in [0,T],\ i=1,\ldots,p,$$
where
$$t^*={\rm argmin}_{t\in[0,T]}\left[\frac 12 \frac{(-\langle \hat  b(t), \xi \rangle+x)^2} {\Big({\sum_{j=1}^p
\xi_j^{\theta/(\theta-2)}k_{jj}(t,t)^{\theta/(\theta-2)}}\Big)^{(\theta-2)/\theta}}\right]^{\theta/(\theta+2)}.$$
We also remark that $k_{ii}(0,0)=0$ for all $i=1,\ldots,p$; so, in both cases, $t^*\in(0,T]$.

\subsection{Exit probabilities from quadrants}\label{sect:LCP-quadrant}
In this section we consider
\begin{equation}\label{eq:LCP-quadrant-def}
p_n:=\P\Big(\bigcap_{i=1}^p\Big\{\sup_{t\in[0,T]} Z^n_i(t) \geq x_i\Big\}\Big),
\end{equation}
for $\hat b_1,\ldots,\hat b_p $ such that
$x_1-\hat b_1(t) ,\ldots,x_p-\hat b_p(t) >0$ for $t\in[0,T]$, and $x_1,\ldots,x_p>0$.
Moreover the set $\cl A$ is defined by
$$\cl A:=\left\{ z \in {C}([0,T], \R^p): \sup_{t\in[0,T]}z_i(t)\geq x_i,\ i=1,\ldots,p\right\}.$$
Then we have
\begin{equation}\label{eq:union-of-paths-quadrant}
\cl A=\bigcup_{t_1,\ldots,t_p\in[0,T]}\cl{A}_{t_1,\ldots,t_p},\ \mbox{where}\
\cl{A}_{t_1,\ldots,t_p}:=\left\lbrace z \in {C}([0,T],\R^p): z_i(t_i)=x_i,\ i=1,\ldots,p\right\rbrace.
\end{equation}

\begin{remark}\label{rem:coincident-infima-quadrant}\rm
Here we check (\ref{eq:coincident-infima}), and we follow the same lines of Remark \ref{rem:coincident-infima-halfspace}
(where there is a different set $\cl A$). Again we have to check an inequality. Moreover we can say again that there 
exists $z^*\in \cl{A}$ such that $I_Z(z^*)=\inf_{z \in \cl{A}}I_Z(z)$ and, to avoid trivialities, we assume that
$z^*\notin \cl{A}^\circ$. In this case we have $\sup_{0\leq t \leq T}z_i^*(t)\geq x_i$ and, moreover, 
$\sup_{0\leq t \leq T}z_i^*(t)=x_i$ if and only if $i$ belongs to a suitable nonempty set of indices $\mathcal{I}$. 
Therefore, for some $t_1^*,\ldots,t_p^*\in[0,T]$, we have
$z_i^*(t_i^*)=x_i>\hat  b_i(t_i^*)$ if $i\in\mathcal{I}$, and $z_i^*(t_i^*)>x_i>\hat  b_i(t_i^*)$ otherwise.
Then, for every $\varepsilon>0$, let $z^{*,\varepsilon}=(z_1^{*,\varepsilon},\ldots,z_p^{*,\varepsilon})$ be defined by
$$z_i^{*,\varepsilon}(t_i):=\left\{\begin{array}{ll}
(1+\varepsilon)(z_i^*(t_i)-\hat  b_i(t_i))+\hat  b_i(t_i)&\ \mbox{if}\ i\in\mathcal{I}\\
z_i^*(t_i)&\ \mbox{otherwise}.
\end{array}\right.$$
So $z^{*,\varepsilon}\in\cl{A}^\circ$. Indeed, if $i\in\mathcal{I}$, we have
\begin{multline*}
\sup_{0\leq t_i \leq T}z_i^{*,\varepsilon}(t_i)
\geq z_i^{*,\varepsilon}(t_i^*)=(1+\varepsilon)(z_i^*(t_i^*)-\hat  b_i(t_i^*))+\hat  b_i(t_i^*)\\
=(1+\varepsilon)(x_i-\hat  b_i(t_i^*))+\hat  b_i(t_i^*)
=x_i+\varepsilon(x_i-\hat  b_i(t_i^*))>x_i;
\end{multline*}
on the other hand, if $i\notin \mathcal{I}$, we have $\sup_{0\leq t_i \leq T}z_i^{*,\varepsilon}(t_i)
\geq z_i^{*,\varepsilon}(t_i^*)=z_i^*(t_i^*)>x_i$. Finally we conclude as we did in Remark \ref{rem:coincident-infima-halfspace}.
\end{remark}

\begin{remark}\label{rem:simultaneous-exit}\rm
One could consider the case in which all the components reach or exceed certain levels simultaneously. In such a case
we have to consider
$$\widehat{p}_n:=\P\Big(\bigcup_{t\in[0,T]}\bigcap_{i=1}^p\Big\{Z^n_i(t) \geq x_i\Big\}\Big)$$
and the set of paths 
$$\widehat{\cl A}=\bigcup_{t_1,\ldots,t_p\in[0,T],t_1=\cdots=t_p}\cl{A}_{t_1,\ldots,t_p}$$
(instead of $p_n$ in (\ref{eq:LCP-quadrant-def}) and $\cl A$ in (\ref{eq:union-of-paths-quadrant})). Obviously, by construction,
we have $\widehat{p}_n\leq p_n$ and $\widehat{\cl A}\subset \cl A$. So we should have
$$\widehat w:=\inf_{z \in \widehat{\cl A}^\circ}I_Z(z)=\inf_{z \in \widehat{\cl A}}I_Z(z),$$
which yields the logarithmic asymptotic estimate
$$\lim_{n\to +\infty}\frac{1}{\gamma_n}\log(\widehat p_n)=-\widehat w.$$
In our opinion the method used in this paper to compute the exponential decay rates do not work well and we are not able to give
a formula for $\widehat w$. However we have the obvious trivial inequality $\widehat w\geq w$.
\end{remark}

In the following Propositions \ref{prop:LCP-quadrant-equal} and \ref{prop:LCP-quadrant-independent} we compute $w$ noting 
that
\begin{equation}\label{eq:w-def-quadrant}
w=\inf_{z \in \cl{A}}I_Z(z)=\inf_{t_1,\ldots,t_p\in[0,T]}\inf_{z \in \cl{A}_{t_1,\ldots,t_p}}I_Z(z);
\end{equation}
moreover, for every $t_1,\ldots,t_p\in [0,T]$, we compute $\inf_{z \in \cl{A}_{t_1,\ldots,t_p}}I_Z(z)$ by taking $z$ in
the dense set $\cl{D}$ defined above (see eq. (\ref{eq:z-dense})) by applying the Lagrange multiplied methods (we 
can do this because the rate function $I_Z$ is strictly convex; see Remark \ref{rem:constraint-on-theta}). Now we are ready to 
prove Proposition \ref{prop:LCP-quadrant-equal} (see Corollary \ref{cor:main-equal} and Remark \ref{rem:main-equal}) and Proposition 
\ref{prop:LCP-quadrant-independent} (see Corollary \ref{cor:main-ind} and Remark \ref{rem:main-ind}). In Proposition 
\ref{prop:LCP-quadrant-independent} we have $\theta_1,\ldots,\theta_p>2$ as we said in Remark \ref{rem:constraint-on-theta} and, 
differently from Proposition \ref{prop:LCP-halfspace-independent}, we do not need to assume that $\theta_1,\ldots,\theta_p$ are all 
coincident.

\begin{proposition}\label{prop:LCP-quadrant-equal}
Let $p_n$ be as in (\ref{eq:LCP-quadrant-def}), where
$$Z^n(t)=\widetilde{A}^nX^n(t)+\hat b(t)\quad(t\in[0,T]),$$
$(X^n)_{n\in\N}$ is a continuous Gaussian process as in Condition (C1), $\hat b$ is as in Remark \ref{rem:main-equal}, and 
$(\widetilde{A}^n)_{n\in\N}$ in Condition (C2) for some $d,\theta>0$. Then (\ref{eq:LCP-limit}) holds with $w=w_Q^{(=)}$, where
$$w_Q^{(=)}:=\psi(d,\theta)\inf_{t_1,\ldots,t_p\in[0,T]}
\left[\frac 12\sum_{i=1}^p \frac{(-\hat b_i(t_i)+x_i)^2}{k_{ii}(t_i,t_i)}\right]^{\theta/(\theta+2)}.$$
\end{proposition}
\begin{proof}
	We have to compute $w$ in (\ref{eq:w-def-quadrant}) where $\cl{A}_{t_1,\ldots,t_p}$ is as in (\ref{eq:union-of-paths-quadrant})
	and $I_Z$ is the rate function $I_Z^{(=)}$ in (\ref{eq:rate-fun-Z-ggG}). Thus, by referring to the elements in the 
	dense set of paths $\cl{D}_{\hat b}$ (see (\ref{eq:z-dense}) and (\ref{eq:z-dense-norm})), we have to minimize
	$$\psi(d,\theta)\left(\sum_{i=1}^p\frac{1}{2}\int_0^T \int_0^T k_{ii}(u,v)\ d\lambda_i(u)\ d\lambda_i(v)\right)^{\theta/(\theta+2)}$$
	with respect to the vector measures $\lambda$, subjected to the following constraint
	$$\int_0^T k_{ii}(t_i,v)\ d\lambda_i(v)=x_i-\hat b_i(t_i),\quad i=1,\ldots,p. $$
	So we use the method of Lagrange multipliers and, for $\gamma_1,\ldots,\gamma_p\in\R$, we have to find the stationary points of
	$$\displaylines{\cl{L}(\lambda,\gamma):=\psi(d,\theta)\left(\sum_{i=1}^p\frac{1}{2}\int_0^T \int_0^T k_{ii}(u,v)
		\ d\lambda_i(u)\ d\lambda_i(v)\right)^{\theta/(\theta+2)}\!\!\!\cr
		-{\sum_{i=1}^p}\gamma_i\Big(\int_0^Tk_{ii}(t_i,v)\  d\lambda_i(v)-x_i+\hat b_i(t_i)\Big).}$$
	Then, for every  $\eta\in  \cl{M}^1[0, 1]$, we have
	$$\displaylines{\psi(d,\theta)\frac{\theta}{\theta+2}\left(\frac{1}{2}\sum_{i=1}^p\int_0^T \int_0^T k_{ii}(u,v)\ d\lambda_i(u)\ d\lambda_i(v)\right)^{-\frac{2}{\theta+2}}\int_0^T\int_0^T k_{ii}(u,v)\ d\lambda_i(u)\ d\eta(v)+
		\cr - \gamma_i\int_0^T k_{ii}(t_i,v)\  d\eta(v)=0,\quad i=1,\ldots,p;}$$
	therefore we have to find a vector measures $\lambda$ such that, for every $v\in [0,T]$,
	$$\displaylines{\psi(d,\theta)\frac{\theta}{\theta+2}\left(\frac{1}{2}\sum_{i=1}^p\int_0^T \int_0^T k_{ii}(u,v)
		\ d\lambda_i(u)\ d\lambda_i(v)\right)^{-\frac{2}{\theta+2}} \int_0^T k_{ii}(u,v)\ d\lambda_i(u)\cr
		- \gamma_i k_{ii}(t_i,v) =0,\quad i=1,\ldots,p.}$$
	A solution is  $\lambda_i=c_i\delta_{\{ t_i\} }$ for some $c_i\in\R$ (for every $i=1,\ldots,p$). So we have to
	consider the following system
	$$\left\{
	\begin{array}{ll}
		\displaystyle \psi(d,\theta)\frac{\theta}{\theta+2}\left(\frac{1}{2}\sum_{i=1}^p c_i^2k_{ii}(t_i,t_i)\right)^{-\frac{2}{\theta+2}} c_i-\gamma_i=0& i=1,\ldots,p\\
		\displaystyle c_i k_{ii}(t_i,t_i)=-\hat b_i(t_i)+x_i.& i=1,\ldots,p.\end{array}
	\right.$$
	In particular we immediately get
	$$\displaylines{c_i=  \frac {-\hat b_i(t_i)+x_i}{k_{ii}(t_i,t_i)},\quad i=1,\ldots,p.}$$
	We have a solution of the Lagrange multipliers problem, which is a critique point for the functional we want to minimize.
	Moreover, since this is a strictly convex functional restricted on a linear subspace of $\cl{M}^p[0, 1]$, it is still strictly
	convex, and therefore the critique point is actually its unique point of minimum. In conclusion we immediately get the desired
	expression for $w$ by (\ref{eq:w-def-quadrant}); indeed we have
	$\inf_{z \in \cl{A}_{t_1,\ldots,t_p}}I_Z(z)= \psi(d,\theta)\left[\sum_{i=1}^p\frac 12 c_i^2 k_{ii}(t_i,t_i)\right]^{\theta/(\theta+2)}$
	with $c_1,\ldots,c_p$ computed above.
\end{proof}

\begin{proposition}\label{prop:LCP-quadrant-independent}
Let $p_n$ be as in (\ref{eq:LCP-quadrant-def}), where
$$Z^n(t)=A^n\circ X^n(t)+\hat b(t)\quad(t\in[0,T]),$$
$(X^n)_{n\in\N}$ is a continuous Gaussian process as in Condition (C1), $\hat b$ is as in Remark \ref{rem:main-ind}, and 
$(A^n)_{n\in\N}=((\widetilde{A}^n_1,\ldots,\widetilde{A}^n_p))_{n \in \N}$, where 
$(\widetilde{A}^n_1)_{n \in \N},\ldots(\widetilde{A}^n_p)_{n \in \N},$ are as in Condition (C2) for some $d_1,\ldots,d_p>0$ and 
$\theta_1,\ldots,\theta_p>2$. Then (\ref{eq:LCP-limit}) holds with $w=w_Q^{(\bot)}$, where
$$w_Q^{(\bot)}:=\inf_{t_1,\ldots,t_p\in[0,T]}\sum_{i=1}^p\psi(d_i,\theta_i)\left[\frac 12 \frac{(-\hat b_i(t_i)+x_i)^2}{k_{ii}(t_i,t_i)}\right]^{\theta_i/(\theta_i+2)}.$$
\end{proposition}
\begin{proof}
We have to compute $w$ in (\ref{eq:w-def-quadrant}) where $\cl{A}_{t_1,\ldots,t_p}$ is as in (\ref{eq:union-of-paths-quadrant})
and $I_Z$ is the rate function $I_Z^{(\bot)}$ in (\ref{eq:rate-fun-Z}). We follow the same lines of the proof of Proposition
\ref{prop:LCP-quadrant-equal} and we omit some details. Then we use again the method of Lagrange multipliers and,
for $\gamma_1,\ldots,\gamma_p\in\R$, we have to find the stationary points of 
$$\displaylines{\cl{L}(\lambda,\gamma):=\sum_{i=1}^p\psi(d_i,\theta_i)\left(\frac{1}{2}\int_0^T \int_0^T k_{ii}(u,v)
	\ d\lambda_i(u)\ d\lambda_i(v)\right)^{\theta_i/(\theta_i+2)}\!\!\!\cr
	-{\sum_{i=1}^p}\gamma_i\Big(\int_0^T k_{ii}(t_i,v)\  d\lambda_i(v)-x_i+\hat b_i(t_i)\Big).}$$
Then, for every  $\eta\in  \cl{M}^1[0, 1]$, we have
$$\displaylines{\psi(d_i,\theta_i)\frac{\theta_i}{\theta_i+2}\left(\frac{1}{2}\int_0^T \int_0^T k_{ii}(u,v)\ d\lambda_i(u)\ d\lambda_i(v)\right)^{-\frac{2}{\theta_i+2}}\int_0^T\int_0^T k_{ii}(u,v)\ d\lambda_i(u)\ d\eta(v)+
	\cr - \gamma_i\int_0^T k_{ii}(t_i,v)\  d\eta(v)=0,\quad i=1,\ldots,p;}$$
therefore we have to find a vector measures $\lambda$ such that, for every $v\in [0,T]$,
$$\displaylines{\psi(d_i,\theta_i)\frac{\theta_i}{\theta_i+2}\left(\frac{1}{2}\int_0^T \int_0^T k_{ii}(u,v)
	\ d\lambda_i(u)\ d\lambda_i(v)\right)^{-\frac{2}{\theta_i+2}} \int_0^T k_{ii}(u,v)\ d\lambda_i(u)\cr
	- \gamma_i k_{ii}(t_i,v) =0,\quad i=1,\ldots,p.}$$
Then, again, a solution is  $\lambda_i=c_i\delta_{\{ t_i\} }$ for some $c_i\in\R$ (for every $i=1,\ldots,p$). So we have to
consider the following system
$$\left\{
\begin{array}{ll}
	\displaystyle \psi(d_i,\theta_i)\frac{\theta_i}{\theta_i+2}\left(\frac{1}{2}k_{ii}(t_i,t_i)\right)^{-\frac{2}{\theta_i+2}} c_i^{(\theta_i-2)/(\theta_i+2)}-\gamma_i=0& i=1,\ldots,p\\
	\displaystyle c_i k_{ii}(t_i,t_i)=-\hat b_i(t_i)+x_i&i=1,\ldots,p.\end{array}
\right.$$
In particular we immediately get
$$\displaylines{c_i=  \frac {-\hat b_i(t_i)+x_i}{k_{ii}(t_i,t_i)},\quad i=1,\ldots,p.}$$
So we can conclude following the same lines of the final part of the proof of Proposition
\ref{prop:LCP-quadrant-equal}. In particular we immediately get the desired expression for
$w$ by (\ref{eq:w-def-quadrant}); indeed we have
$\inf_{z \in \cl{A}_{t_1,\ldots,t_p}}I_Z(z)=\sum_{i=1}^p \psi(d_i,\theta_i)\left[\frac 12 c_i^2 k_{ii}(t_i,t_i)\right]^{\theta_i/(\theta_i+2)}$
with $c_1,\ldots,c_p$ computed above.
\end{proof}

\paragraph{On the most likely paths leading to exit.}
This is the analogue of the discussion presented in Section \ref{sect:LCP-halfspace}. A closer look at of the proofs of Propositions
\ref{prop:LCP-quadrant-equal} and \ref{prop:LCP-quadrant-independent} reveal that, in both cases, we can define a most likely path 
$z^*$ as follows
$$z_i^*(u):=\frac{-\hat b_i(t_i^*)+x_i}{k_{ii}(t_i^*,t_i^*)}k_{ii}(u,t_i^*)+\hat b_i(u), \quad u \in [0,T],\ i=1,\ldots,p,$$
with two different definitions of $(t^*_1,\ldots,t^*_p)$. In the case of Proposition \ref{prop:LCP-quadrant-equal} we have
$$(t^*_1,\ldots,t^*_p)={\rm argmin}_{t_1,\ldots,t_p\in[0,T]}\left[\frac 12\sum_{i=1}^p
\frac{(-\hat b_i(t_i)+x_i)^2}{k_{ii}(t_i,t_i)}\right]^{\theta/(\theta+2)},$$
and in the case of Proposition \ref{prop:LCP-quadrant-independent} we have
$$(t^*_1,\ldots,t^*_p)={\rm argmin}_{t_1,\ldots,t_p\in[0,T]}
\sum_{i=1}^p\psi(d_i,\theta_i)\left[\frac 12 \frac{(-\hat b_i(t_i)+x_i)^2}{k_{ii}(t_i,t_i)}\right]^{\theta_i/{(\theta_i+2)}},$$
that is
$t^*_i={\rm argmin}_{t_i\in[0,T]}\left[\frac 12 \frac{(-\hat b_i(t_i)+x_i)^2}{k_{ii}(t_i,t_i)}\right]^{\theta_i/{(\theta_i+2)}}$
for every $i=1,\ldots,p$. We also remark that $k_{ii}(0,0)=0$ for all $i=1,\ldots,p$; so, in both cases, $t_i^*\in(0,T]$
for all $i=1,\ldots,p$.

\section{Comparisons between asymptotic rates}\label{sect:final}
Inequalities between rate functions allow to compare the convergence of two sequences of stochastic processes. Moreover, by 
taking into account the limit (\ref{eq:LCP-limit}), inequalities between rate functions allow to compare the exponential 
decay rates of exit probabilities (indeed, as we pointed out, $w$ is equal to the infimum of the rate function over a suitable 
set of paths). In view of the results presented in this section it is useful to refer to the following lemma.

\begin{lemma}\label{lemma:inequalities}
Let $x_1,\ldots, x_n$ be non negative number. Then:
for $r\in (0,1)$ $(r\in (1,+\infty),\, resp.)$
$$ \sum_{i=1}^p x_i^r\geq \Big(\sum_{i=1}^p x_i\Big)^r \quad \Big(\sum_{i=1}^p x_i^r\leq \Big(\sum_{i=1}^p x_i\Big)^r\, resp.\Big),$$
and the equality holds if and only if the set $\{i\in\{1,\ldots,p\}:x_i\neq 0\}$
has at most one element.
\end{lemma}
\begin{proof}
The case $x_1=\ldots=x_n=0$ is trivial.
We remark that $y^r\geq y$ for $y\in[0,1]$ if $r\in(0,1)$, and $y^r\leq y$ for $y\in[0,1]$ if 
$r\in(1,\infty)$; moreover we have $y^r=y$ if and only if $y=0$ or $y=1$. Now we assume that $x_i\neq 0$
for some $i\in\{1,\ldots,p\}$, and we have ${\cal S}(x)=\sum_{i=1}^p x_i>0$; then, if $r\in(0,1)$, we have
$$\sum_{i=1}^p x_i^r = {\cal S}(x)^r\sum_{i=1}^p\Big(\frac {x_i} {{\cal S}(x)}\Big)^r
\geq{\cal S}(x)^r\sum_{i=1}^p\frac {x_i} {{\cal S}(x)}= {\cal S}(x)^r,$$
and the inverse inequality for $r\in(1,+\infty)$. Finally it is easy to see that the inequalities
turn into equalities if and only if the set $\{i\in\{1,\ldots,p\}:x_i\neq 0\}$ has at most one element.
\end{proof}

The results in this section are collected in the following proposition. We remark that, since we want to compare rate functions and 
exponential decay rates concerning different situations, we always have $d_1=\cdots=d_p=d$ and $\theta_1=\cdots=\theta_p=\theta$ for
some $d>0$ and $\theta>2$ (actually, for the statement (ii), this is a consequence of the hypotheses in Proposition 
\ref{prop:LCP-halfspace-independent}).

\begin{proposition}\label{prop:inequalities}
The following three statements hold.\\
(i) Let $I_Z^{(=)}$ be the rate function in (\ref{eq:rate-fun-Z-ggG}), and let $I_Z^{(\bot)}$ be the rate function in (\ref{eq:rate-fun-Z})
with $d_1=\cdots=d_p=d$ and $\theta_1=\cdots=\theta_p=\theta$ for some $d>0$ and $\theta>2$. Then $I_Z^{(=)}(z)\leq I_Z^{(\bot)}(z)$ for all
$z\in C([0,T],\R^p)$. Moreover the equality holds if and only if the set $\{i\in\{1,\ldots,p\}:z_i\neq \hat b_i\}$ has at most one element.\\
(ii) Let $w_H^{(=)}$ and $w_H^{(\bot)}$ be the exponential decay rates in Propositions \ref{prop:LCP-halfspace-equal} and 
\ref{prop:LCP-halfspace-independent}. Then $w_H^{(=)}\leq w_H^{(\bot)}$. Moreover the equality holds if and only if the set 
$\{i\in\{1,\ldots,p\}:\xi_i\neq 0\}$ has at most one element.\\
(iii) Let $w_Q^{(=)}$ and $w_Q^{(\bot)}$ be the exponential decay rates in Propositions \ref{prop:LCP-quadrant-equal} and 
\ref{prop:LCP-quadrant-independent}, and take the second one with $d_1=\cdots=d_p=d$ and $\theta_1=\cdots=\theta_p=\theta$ for some $d>0$ and 
$\theta>2$. Then $w_Q^{(=)}\leq w_Q^{(\bot)}$. Moreover the equality holds if and only if $p=1$.
\end{proposition}
\begin{proof}
We prove the statements separately.\\
(i) We trivially have $I_Z^{(=)}(z)=I_Z^{(\bot)}(z)=0$ if $z=\hat b$, i.e. $z=(\hat b_1,\ldots,\hat b_p)$. If 
$z\neq\hat b$, we can apply Lemma \ref{lemma:inequalities} with $x_i=\frac 12 \lVert z_i-\hat b_i\rVert^{2}_{{\cl{H}_{i}}}$ 
for $i\in\{1,\ldots,p\}$ and $r=\frac{\theta}{\theta+2}\in(0,1)$.\\
(ii) By Lemma \ref{lemma:inequalities} with $x_i= \xi_i^{2}k_{ii}(t,t)$ for $i\in\{1,\ldots,p\}$ and 
$r=\frac{\theta}{\theta-2}\in(1,+\infty)$ we get
$$\left(\sum_{i=1}^p \xi_i^{2}k_{ii}(t,t)\right)^{\theta/(\theta-2)}
\geq\sum_{i=1}^p\xi_i^{2\theta/(\theta-2)}k_{ii}(t,t)^{\theta/(\theta-2)};$$
then we obtain
$$\left[\frac 12 \frac{(-\langle \hat  b(t), \xi \rangle+x)^2} {\sum_{j=1}^p \xi_j^{2}k_{jj}(t,t)}\right]^{\theta/(\theta+2)}
\leq\left[\frac 12 \frac{(-\langle \hat  b(t), \xi \rangle+x)^2} {\Big(\sum_{j=1}^p \xi_j^{2\theta/(\theta-2)}
k_{jj}(t,t)^{\theta/(\theta-2)}\Big)^{(\theta-2)/\theta}}\right]^{\theta/(\theta+2)},$$
and the equality holds if and only if the set $\{i\in\{1,\ldots,p\}:\xi_i\neq 0\}$
has one element. If $\{i\in\{1,\ldots,p\}:\xi_i\neq 0\}$ has at least two elements, then we get the desired
strict inequalities between the infima with respect to $t\in[0,T]$ (in both left and right hand sides the 
infima are attained).\\
(iii) By Lemma \ref{lemma:inequalities} with $x_i=\frac 12 \frac{(-\hat b_i(t_i)+x_i)^2}{k_{ii}(t_i,t_i)}$
for $i\in\{1,\ldots,p\}$ (note that they are all positive) and $r=\frac{\theta}{\theta+2}\in(0,1)$ we get
$$\displaylines
{\sum_{i=1}^p\left[\frac 12 \frac{(-\hat b_i(t_i)+x_i)^2}{k_{ii}(t_i,t_i)}\right]^{\theta/(\theta+2)}
\geq
\left[\sum_{i=1}^p\frac 12 \frac{(-\hat b_i(t_i)+x_i)^2}{k_{ii}(t_i,t_i)}\right]^{\theta/(\theta+2)},}$$
and the equality holds if and only if $p=1$. If $p\geq 2$, then we get the desired strict inequalities 
between the infima with respect to $t_1,\ldots,t_p\in[0,T]$ (in both left and right hand sides the infima are attained).
\end{proof}

The inequalities between exponential decay rates of exit probabilities (see statements (ii) and (iii) in
Proposition \ref{prop:inequalities}) have some analogies with some inequalities in the literature. Here we recall
the inequalities in Proposition 4.1 in \cite{MST} concerning a multivariate risk process with delayed claims in
insurance; in such a case the exponential decay rate is called Lundberg parameter. In particular one has a smaller
exponential decay rate when the joint distribution of the claims is larger with respect to the supermodular order.
We also recall that, as discussed in \cite{MST} (see Section 4.2), comonoticity is the strongest case of dependence
with respect to the supermodular order (see \cite{MullerStoyan}, Theorem 3.9.5 (c), p. 114; see also
\cite{BauerleMuller}, Theorem 2.1) when one compares two joint distributions with the same marginal distributions; 
for this reason in Proposition \ref{prop:inequalities}(iii) we assume that $d_1=\cdots=d_p=d$ and 
$\theta_1=\cdots=\theta_p=\theta$. We also recall that $Z^n=A^n X^n +B^n$ in Theorem \ref{th:main-equal} can be
rewritten as $Z^n=(A^n,\ldots,A^n)\circ X^n +B^n$ and the random vector $(A^n,\ldots,A^n)$ is comonotonic because 
the components are all coincident.

\section{A discussion on Condition (C2)}\label{sect:discussion-C2}
In this section we present Proposition \ref{prop:sufficient-condition-for-C2}, which yields Condition (C2). Moreover
we present two examples; in particular Example \ref{ex:L-beta-continuation} with $\rho=1/2$ concerns 
the \emph{vector valued generalized grey Brownian motion} (vggBM) in \cite{BockGrothausOrge} and another
$p$-variate process that, for $p=1$, reduces to the \emph{generalized grey Brownian motion} (ggBM) in \cite{BockDesmettreDasilva}.
We recall that $(\sqrt{L_\beta}X(t))_{t\in[0,T]}$ is a ggBM if $(X(t))_{t\in[0,T]}$ is a fractional Brownian motion 
with Hurst parameter $\frac{\alpha}{2}$, with $\alpha\in(0,2)$, indipendent of the random variable $L_\beta$ as in Example
\ref{ex:L-beta} (see below). We also recall that $(X(t))_{t\in[0,T]}$ is a fractional Brownian motion with Hurst parameter 
$\frac{\alpha}{2}$, with $\alpha\in(0,2)$, if we have the following covariance functions
$$k_{ij}(t,s)=\left\{\begin{array}{ll}
	\frac 12(t^\alpha+s^\alpha-|t-s|^\alpha)&\ \mbox{if}\ i=j\\
	0&\ \mbox{otherwise}.
\end{array}\right.$$
A reference for the ggBM is Proposition 3 in \cite{MuraPagnini}.

\begin{proposition}\label{prop:sufficient-condition-for-C2}
Let $\widetilde{A}$ be a continuous positive random variables such that
$$P(\widetilde{A}\geq r)=e^{-dr^\theta(1+o(1))}$$
for some $d,\theta>0$. Moreover set $\widetilde{A}_n:=\gamma_n^{-1/\theta}\widetilde{A}$. Then Condition (C2) holds.
\end{proposition}

\begin{proof}
	We start with the upper bound for closed sets
	$$\limsup_{n \rightarrow +\infty}\frac{1}{\gamma_n}\log P(\gamma_n^{-1/\theta}\widetilde{A}\in C)\leq-\inf_{x\in C}{\cal J}(x)$$
	for every closed subset $C$ of $[0,+\infty)$. This immediately holds if $C$ is empty or if $0\in C$. On the contrary, if $C$
	is nonempty and $0\notin C$, there exists $x_C>0$ such that $x_C:=\inf C$ and $x_C\in C$ and we get
	$$P(\gamma_n^{-1/\theta}\widetilde{A}\in C)\leq P(\gamma_n^{-1/\theta}\widetilde{A}\geq x_C)
	=P(\widetilde{A}\geq \gamma_n^{1/\theta}x_C),$$
	which yields
	$$\displaylines{\limsup_{n \rightarrow +\infty}\frac{1}{\gamma_n}\log P(\gamma_n^{-1/\theta}\widetilde{A}\in C)\leq
		\limsup_{n \rightarrow +\infty}\frac{1}{\gamma_n}\log P(\widetilde{A}\geq \gamma_n^{1/\theta}x_C)\cr
		=\limsup_{n \rightarrow +\infty}\frac{-d(\gamma_n^{1/\theta}x_C)^\theta(1+o(1))}{\gamma_n}=-dx_C^\theta
		=-{\cal J}(x_C)=-\inf_{x\in C}{\cal J}(x)}$$
	(in the last equality we take into account that ${\cal J}$ is continuous and increasing).
	
	For the lower bound for the open sets it is well-known that it suffices to show that, for every $x\geq 0$ and for every open 
	set $O$ such that $x\in O$, we have
	$$\liminf_{n \rightarrow +\infty}\frac{1}{\gamma_n}\log P(\gamma_n^{-1/\theta}\widetilde{A}\in O)\geq-{\cal J}(x).$$
	For $x=0$ this estimate holds because $\gamma_n^{-1/\theta}\widetilde{A}$ tends to zero (almost surely). For $x>0$, 
	since we can find $\varepsilon>0$ such that $(x-\varepsilon,x+\varepsilon)\subset O\cap(0,\infty)$, we have
	$$P(\gamma_n^{-1/\theta}\widetilde{A}\in O)\geq P(\gamma_n^{-1/\theta}\widetilde{A}\in (x-\varepsilon,x+\varepsilon))
	=P(\gamma_n^{-1/\theta}\widetilde{A}\geq x-\varepsilon)-P(\gamma_n^{-1/\theta}\widetilde{A}\geq x+\varepsilon);$$
	moreover, since
	$$\liminf_{n \rightarrow +\infty}\frac{1}{\gamma_n}\log(P(\gamma_n^{-1/\theta}\widetilde{A}\geq x-\varepsilon)
	-P(\gamma_n^{-1/\theta}\widetilde{A}\geq x+\varepsilon))\geq-d(x-\varepsilon)^\theta$$
	by Lemma 19 in \cite{GaneshTorrisi}, we get
	$$\liminf_{n \rightarrow +\infty}\frac{1}{\gamma_n}\log P(\gamma_n^{-1/\theta}\widetilde{A}\in O)
	\geq-d(x-\varepsilon)^\theta,$$
	and we obtain the desired estimate by letting $\varepsilon$ go to zero.
\end{proof}

Now we present an example (see Proposition 3 and eq. (A.1) in \cite{MuraPagnini}; see also eq. (A.4) in \cite{MuraPagnini}
for the moment generating function).

\begin{example}\label{ex:L-beta}\rm
	For $\beta\in(0,1)$, let $L_\beta$ be a random variable with density
	$$M_\beta(\tau)=\sum_{k=0}^{\infty}\frac{(-\tau)^k}{k!\Gamma(-\beta k+1-\beta)}.$$
	Then we have
	\begin{equation}\label{eq:L-beta-mgf}
		\mathbb{E}\left[\exp(\eta L_\beta)\right]=E_\beta(\eta)\quad (\mbox{for all}\ \eta\in\R),
	\end{equation}
	where $E_\beta(z):=\sum_{h=0}^\infty\frac{z^h}{\Gamma(\beta h+1)}$ is the Mittag-Leffler function. Note that 
	equation (A.4) in \cite{MuraPagnini} and some other references provides this formula only for $\eta\leq 0$; however
	this restriction is not needed because we can refer to the analytic continuation of the Laplace transform with 
	complex argument. 
\end{example}

The next Proposition \ref{prop:L-beta} provides a LDP as the one stated in Proposition \ref{prop:sufficient-condition-for-C2}
(with $\widetilde{A}=L_\beta$ as in Example \ref{ex:L-beta}).

\begin{proposition}\label{prop:L-beta}
	Let $L_\beta$ be a random variable as in Example \ref{ex:L-beta}, for some $\beta\in(0,1)$.
	Then $(\gamma_n^{-1+\beta}L_\beta)_{n\in \N}$ satisfies the LDP, on $[0,+\infty)$, with the speed $\gamma_n$ and good rate 
	function $\cal J$ defined by (\ref{eq:rf-C2}) in Condition (C2), with $d:=\beta^{\beta/(1-\beta)}-\beta^{1/(1-\beta)}$
    and $\theta:=\frac{1}{1-\beta}$.
\end{proposition}
\begin{proof}
	The proof consists of a straightforward application of the G\"artner Ellis Theorem (see e.g. Theorem 2.3.6 in
	\cite{DemboZeitouni}). Then, by taking into account the moment generating function
	in (\ref{eq:L-beta-mgf}) and the asymptotic behavior of the Mittag-Leffler function (see e.g. (3.4.14)-(3.4.15) 
	in Proposition 3.6 in \cite{GorenfloKilbasMainardiRogosin}), we have
	$$\Lambda(\eta)=\lim_{n\to +\infty}\frac{1}{\gamma_n}\log\mathbb{E}\left[\exp((\gamma_n\,\gamma_n^{-1+\beta}\eta L_\beta))\right]
	=\lim_{n\to +\infty}\frac{1}{\gamma_n}\log E_\beta(\gamma_n^\beta\eta)
	=\left\{\begin{array}{ll}
			\eta^{1/\beta}&\ \mbox{if}\ \eta\geq 0\\
			0&\ \mbox{otherwise}.
	\end{array}\right.$$
	Then, since the function $\Lambda$ is differentiable, the desired LDP holds with good rate function $\Lambda^*$ defined by
	$$\Lambda^*(x):=\sup_{\eta\in\mathbb{R}}\{\eta\, x-\Lambda(\eta)\},$$
	which coincides with the rate function ${\cal J}$ in the statement of the proposition (with 
	$d:=\beta^{\beta/(1-\beta)}-\beta^{1/(1-\beta)}$ and $\theta:=\frac{1}{1-\beta}$); indeed this can be readily checked
	(note that, actually, we should consider $x\in\mathbb{R}$; however here we already know that we can neglect the case $x<0$ 
	because we deal with non-negative random variables).
\end{proof}

\begin{example}\label{ex:L-beta-continuation}\rm
	One could try to consider $\widetilde{A}=(L_\beta)^\rho$ for some $\rho>0$ in place of $\widetilde{A}=L_\beta$.
	Then, if we combine the LDP in Proposition \ref{prop:L-beta} with the contraction principle (see e.g. Theorem 4.2.1 in 
	\cite{DemboZeitouni}), we can easily check that $((\gamma_n^{-1+\beta}L_\beta)^\rho)_{n\in \N}$ satisfies the LDP, on 
	$[0,+\infty)$, with the speed $\gamma_n$ and good rate function $\cal J$ defined by (\ref{eq:rf-C2}) in Condition (C2), 
	with $d:=\beta^{\beta/(1-\beta)}-\beta^{1/(1-\beta)}$ and $\theta:=\frac{1}{\rho(1-\beta)}$. So, when $\rho=1/2$,
	we can apply apply Corollary \ref{cor:main-ind} to the vggBM in \cite{BockGrothausOrge} and Corollary 
	\ref{cor:main-equal} to some processes that reduces to the ggBM in \cite{BockDesmettreDasilva} when $p=1$.
\end{example}

\paragraph{Funding.}
The authors are supported by MIUR Excellence Department Project awarded to the Department of Mathematics, University
of Rome Tor Vergata (CUP E83C18000100006 and CUP E83C23000330006), by University of Rome Tor Vergata (project
"Asymptotic Methods in Probability" (CUP E89C20000680005) and project "Asymptotic Properties in Probability" \\
(CUP E83C22001780005)) and by Indam-GNAMPA.

\paragraph{Acknowledgements.}
The authors thank two referees for the careful reading of an earlier version of the paper and for some useful comments.
Moreover the authors also thank Francesco Mainardi, Gianni Pagnini and Vittoria Sposini for some hints and comments 
on the literature about the generalized grey Brownian motion.


\end{document}